\documentclass[12pt]{article}
\usepackage{amsmath,amsthm,amsfonts,amscd}
\title {Rationality of moduli of vector bundles on curves}
\author {Aidan Scho{f}\i eld and Alastair King}
\date {}
\numberwithin{equation}{section}
\theoremstyle{plain}
\newtheorem{theorem}{Theorem}[section]
\newtheorem{lemma}[theorem]{Lemma}

\newtheorem{proposition}[theorem]{Proposition}
\theoremstyle{definition}
\newtheorem{definition}[theorem]{Definition}

\DeclareMathAlphabet{\mathsc}{OT1}{cmr}{m}{sc}
\DeclareMathOperator{\rk}{rk}
\DeclareMathOperator{\Hom}{Hom}

\DeclareMathOperator{\Ext}{Ext}

\DeclareMathOperator{\ext}{ext}
\DeclareMathOperator{\im}{im}
\DeclareMathOperator{\coker}{coker}
\DeclareMathOperator{\Aut}{Aut}

\newcommand{\hcf}[2]{\operatorname{hcf}(#1,#2)}
\newcommand{\ses}[3]{0 \to #1 \to #2 \to #3 \to 0}
\newcommand{\dual}{^{\vee}}
\newcommand{\dsh}{^{\prime}}
\newcommand{\ramto}{\dashrightarrow}
\newcommand{\thetahat}{{\widehat\theta}}
\newcommand{\muhat}{{\widehat\mu}}

\newcommand{\PP}{{\mathbb P}}
\newcommand{\ZZ}{{\mathbb Z}}
\newcommand{\isom}{\cong}
\newcommand{\moduli}[2]{\mathfrak{M}_{#1,#2}}
\newcommand{\parmod}[3]{\mathfrak{P}_{#1,#2,#3}}
\newcommand{\pullback}{\widehat{\mathfrak{P}}}
\newcommand{\quot}[3]{\mathsc{Quot}(#1,#2,#3)}
\newcommand{\Gr}[2]{\mathsc{Gr}({#1},{#2})}
\newcommand{\PGL}[1]{PGL_{#1}}
\newcommand{\GL}[1]{GL_{#1}}
\newcommand{\eulchar}[2]{\chi({#1},{#2})}
\newcommand{\eval}[2]{\varepsilon_{#1}({#2})}
\newcommand{\funfield}[1]{k(#1)}
\newcommand{\Br}[3]{{\mathfrak{Br}}\left(#1\right)}
\newcommand{\Brr}[1]{{\mathfrak{Br}}\left(#1\right)}
\newcommand{\brauer}[2]{\psi_{#1,#2}}
\newcommand{\Fr}[2]{{\mathsc{Fr}}(#1,#2)}
\newcommand{\Frdual}[2]{{\mathsc{Fr}}\dual(#1,#2)}
\newcommand{\cO}{{\mathcal{O}}}
\makeatletter
\newlength{\typesize}
\makeatother

\newlength{\vvoff}
\newlength{\hhoff}

\newcommand{\locateoffcenter}[1]{%
\addtolength{\vvoff}{-0.25\typesize}%
\raisebox{\vvoff}{\hspace{\hhoff}\makebox(0,0){\smash{#1}}}
}
\newcommand{\object}[1]{%
\setlength{\vvoff}{0pt}%
\setlength{\hhoff}{0pt}%
\locateoffcenter{#1}
}
\newcommand{\rlabel}[1]{%
\setlength{\vvoff}{0.75\typesize}%
\setlength{\hhoff}{0pt}%
\locateoffcenter{#1}
}
\newcommand{\dlabel}[1]{%
\setlength{\vvoff}{0pt}%
\setlength{\hhoff}{0.75\typesize}%
\locateoffcenter{#1}
}
\newcommand{\nelabel}[1]{%
\setlength{\vvoff}{0.5\typesize}%
\setlength{\hhoff}{-0.75\typesize}%
\locateoffcenter{#1}
}
\newcommand{\nwlabel}[1]{%
\setlength{\vvoff}{-0.5\typesize}%
\setlength{\hhoff}{-0.75\typesize}%
\locateoffcenter{#1}
}
\newcommand{\swlabel}[1]{%
\setlength{\vvoff}{-0.5\typesize}%
\setlength{\hhoff}{0.75\typesize}%
\locateoffcenter{#1}
}
\begin{document} \maketitle
\section{Introduction} \label{sec:intro}

Let $C$ be a smooth projective curve of genus $g$ over an
algebraically closed field $k$. 
Let $\moduli{r}{d}$ be the moduli space of
stable vector bundles of rank $r$ and degree $d$ over $C$.
This is a smooth quasi-projective variety of dimension $r^2(g-1)+1$,
which is projective when $r$ and $d$ are coprime.
Up to isomorphism, it depends only on the congruence class
of $d$ mod $r$. 
The rank 1 case $\moduli{1}{d}$ is isomorphic to the Jacobian $J(C)$ and 
every moduli space comes equipped with a determinant map
$\det\colon\moduli{r}{d}\to\moduli{1}{d}$ whose fibre
over $L$ is $\moduli{r}{L}$, the moduli space of bundles with fixed
determinant $L$.

The goal of this paper is to describe these moduli spaces in the
birational category, that is, to describe their function fields.
We shall prove the following result.
\begin{theorem}\label{thm:1}
  The moduli space $\moduli{r}{d}$ is birational to 
  $\moduli{h}{0}\times \mathbb{A}^{(r^{2}-h^{2})(g-1)}$,
  where $h=\hcf{r}{d}$. 
\end{theorem}
In other words, there is a dominant rational map
$\mu\colon \moduli{r}{d}\ramto\moduli{h}{0}$
whose generic fibre is rational.
We shall observe that this map restricts to a map between
fixed determinant moduli spaces (not necessarily with the same
determinant) and so, in the case when $r$ and $d$ are coprime,
we obtain the following long believed corollary,
which has been proved in special cases 
(\cite{New75},\cite{New80},\cite{BodYok}).
\begin{theorem}\label{thm:2}
  If $L$ is a line bundle of degree $d$ coprime to $r$,
  then $\moduli{r}{L}$ is a rational variety.
\end{theorem}

To ease the discussion, we use the following terminology
 to describe the relationship between $\moduli{r}{d}$ and $\moduli{h}{0}$.
An irreducible algebraic variety $X$ is \emph{birationally linear}
over another irreducible algebraic variety $Y$ if there exists
a dominant rational map $\phi\colon X\ramto Y$
whose generic fibre is rational, that is, the function field
$\funfield{X}$ is purely transcendental over the function field $\funfield{Y}$.
Such a map $\phi$ will also be called \emph{birationally linear}.

What we shall actually prove is a stronger statement that
the map $\mu$ is birationally linear 
\emph{and} preserves a suitable Brauer class. 
More precisely, for each type $(r,d)$ with $\hcf{r}{d}=h$, 
the moduli space $\moduli{r}{d}$ carries a Brauer class
$\brauer{r}{d}$ for its function field, represented by a central
simple algebra of dimension $h^2$,\
and the map $\mu$ has the property that
$\mu^*(\brauer{h}{0})=\brauer{r}{d}$.
This strengthening of the statement is the key to the proof, 
because it enables an induction on the rank $r$.

In section \ref{sec:step1}, we construct an open dense subvariety of
$\moduli{r}{d}$ as a quotient space of a suitable variety $X_{r,d}$ by
a generically free action of $\PGL{h}$ where $h=\hcf{r}{d}$. This
arises because we are able to show that a general vector bundle $E$ of
type $(r,d)$ arises as a quotient of a particular bundle $F^{h}$ in a
unique way so that we can take $X_{r,d}$ to be a suitable open
subvariety of $\quot{F^h}{r}{d}$ on which $\PGL{h}$ acts in the natural
way. We also arrange in this section that the kernel of the surjection
from $F^{h}$ to $E$ should have smaller rank (at least in the case
where $h\neq r$) and this induces a rational map from $\moduli{r}{d}$
to $\moduli{r_{1}}{d_{1}}$ for some type $(r_{1},d_{1})$ where
$r_{1}<r$. After this, in section \ref{s2}, we show how this
description of $\moduli{r}{d}$ as a quotient space for a generically
free action of $\PGL{h}$ allows us to associate a Brauer class to its
function field and we use this Brauer class to describe birationally
the rational map from $\moduli{r}{d}$ to $\moduli{r_{1}}{d_{1}}$ in
terms of ``twisted Grassmannian varieties''.  In section
\ref{sec:Hecke}, we use parabolic moduli spaces which give us other 
``twisted Grassmannian varieties'' which we may choose to be twisted
``in the same way'' as our map between moduli spaces above. In the
final section, we put these various results together to construct a
birationally linear rational map from $\moduli{r}{d}$ to
$\moduli{h}{0}$. 

\section{The first step} \label{sec:step1}

The purpose of this section is to show that the general bundle $E$
of rank $r$ and degree $d$ may be constructed as a quotient of $F^h$,
where $F$ is a \emph{fixed} bundle of an appropriate type
and $h=\hcf{r}{d}$.

This will enable us to define the Brauer class on $\moduli{r}{d}$
that will be the focus of most of the paper.
Furthermore, we will see that the kernel of the quotient map $q:F^h\to E$
is also general so that we may define a dominant rational map
from $\moduli{r}{d}$ to $\moduli{r_1}{d_1}$,
where $r_1<r$ when $r$ does not divide $d$.
This will be the basis of the inductive construction 
of the birationally linear map to $\moduli{h}{0}$.

We will say that a vector bundle $E$ of rank $r$ and degree $d$ has 
`type' $(r,d)$.
For $E$ of type $\alpha=(r_E,d_E)$
and $F$ of type $\beta=(r_F,d_F)$, we write
\begin{eqnarray*}
  \eulchar{F}{E}&=&\hom(F,E)-\ext(F,E)\\
  &=& r_Fd_E-r_Ed_F-r_Er_F(g-1) = \eulchar{\beta}{\alpha},
\end{eqnarray*}
where $\hom(F,E)=\dim\Hom(F,E)$ and $\ext(F,E)=\dim\Ext(F,E)$.
The middle equality is the Riemann-Roch Theorem.

We start with a lemma about the nature of generic maps between
generic vector bundles.
The proof is closely based on the proof given by Russo \& Teixidor
(\cite{RusTei97} Theorem 1.2) that the tensor product of 
generic bundles is not special; a result originally due to Hirschowitz. 

\begin{lemma}\label{lem:hirsch}
Let $E,F$ be generic vector bundles of fixed types.
Suppose that there exists a non-zero map $\phi\colon F\to E$
and take $\phi$ to be a generic such map.
Then $\ext(F,E)=0$ and $\phi$ has maximal rank.
If $r_E \neq r_F$, then $\coker\phi$ is torsion-free;
in particular, if $r_E < r_F$, then $\phi$ is surjective and if
$r_{E}>r_{F}$ then $\phi$ is injective.
\end{lemma}

\begin{proof}
Let $[\phi]$ denote the homothety class of $\phi$
in $\PP(\Hom(F,E))$.
Then the triple $(F,E,[\phi])$ depends on
$$ p_0= 1-\eulchar{F}{F}+1-\eulchar{E}{E}+\hom(F,E)-1 $$
parameters (cf. \cite{BraGrzNew97} Section 4).
Let $I=\im\phi$, $K=\ker\phi$, $Q=\coker\phi$
and $T$ be the torsion subsheaf of $Q$.
Further, let $Q'=Q/T$ and $I'$ be the inverse image of $T$ in $E$.
Thus we have three short exact sequences
\begin{gather}
0\to K\to F\to I\to 0 \label{eq:KFI}\\
0\to I\to E\to Q\to 0 \label{eq:IEQ}\\
0\to I'\to E\to Q'\to 0 \label{eq:I'EQ'}
\end{gather}
in which all terms except $Q$ are vector bundles.
The triple $(F,E,[\phi])$ is determined by the first and last sequences
(up to homothety) and a map $t\colon I\to I'$ whose cokernel is $T$.
The triple $(I,I',[t])$ depends on $1-\eulchar{I}{I}+r_Id_T$ parameters
and so the whole configuration depends on at most
\begin{gather*}
  p_1= 1-\eulchar{K}{K}+1-\eulchar{I}{I}+r_Id_T +1-\eulchar{Q'}{Q'} \\
  +\ext(I,K)-1+\ext(Q',I')-1
\end{gather*}
parameters.
Now, $E$ and $F$ are stable,
so $\hom(I,K)=\hom(Q',I')=0$ and hence
$\ext(I,K)=-\eulchar{I}{K}$ and $\ext(Q',I')=-\eulchar{Q'}{I'}$.
Furthermore $\eulchar{Q'}{Q'}=\eulchar{Q}{Q}$ and
\begin{equation*}
 \eulchar{Q'}{I'}=r_Q(d_I+d_T)-r_I(d_Q-d_T)-r_Qr_I(g-1)
 =\eulchar{Q}{I} +r_Ed_T
\end{equation*}
Hence, using the bilinearity of $\chi$ in short exact sequences, we get
$$ p_1=1-\eulchar{F}{K}-\eulchar{Q}{E}-\eulchar{I}{I}-r_Qd_T $$
But now $p_0\leq p_1$ and so
\begin{eqnarray*}
 \hom(F,E) &\leq& \eulchar{F}{I}+\eulchar{I}{E}-\eulchar{I}{I} -r_Qd_T\\
           &=& \eulchar{F}{E} - \eulchar{K}{Q} -r_Qd_T
\end{eqnarray*}
and hence 
\begin{equation}\label{eq:extFE}
\ext(F,E)\leq -\eulchar{K}{Q} -r_Qd_T
\end{equation}
A simple dimension count (cf. \cite{Lan83} Lemma 2.1) shows that,
 for general $E$ and $F$ to appear in sequences (\ref{eq:KFI}) and
(\ref{eq:I'EQ'}), it is necessary that
$\eulchar{K}{I}\geq0$ and $\eulchar{I'}{Q'}\geq0$.
In other words,
\begin{gather*}
 r_Kd_I-r_Id_K\geq r_Kr_I(g-1)\\
r_Id_Q-r_Qd_I \geq r_Ir_Q(g-1) + r_Ed_T
\end{gather*}
and hence
\begin{eqnarray*}
\eulchar{K}{Q} &=&
\bigl( r_K(r_Id_Q-r_Qd_I)+ r_Q(r_Kd_I-r_Id_K)- r_Kr_Ir_Q(g-1) \bigr)/r_I\\
&\geq& r_Kr_Q(g-1) +d_T (r_Kr_E/r_I)
\end{eqnarray*}
Thus, subsituting this into (\ref{eq:extFE}), we finally deduce that
$$ \ext(F,E)\leq - r_Kr_Q(g-1) - d_T(r_Q+r_Kr_E/r_I) $$
This is only possible if
(i) $\ext(F,E)=0$, (ii) $r_K=0$ or $r_Q=0$ and 
(iii) unless $r_Q=r_K=0$, we also have $d_T=0$.
But (ii) means that $r_I$ has maximal rank and then
(iii) means that $\coker\phi$ is torsion-free, unless $r_E=r_F$.
\end{proof}

We shall also use the following lemma which may be thought
of as a generalisation of the result that any (bounded) family
of bundles on a curve may be extended to an irreducible family
(cf. \cite{NarRam75} Proposition 2.6).

\begin{lemma}
\label{lem:genext}
Let $\{\mathcal{G}_{x}:x\in X\}$ be an irreducible family of vector
bundles over $C$ and let $\{\mathcal{E}_{y}:y\in Y\}$ be any family
of vector bundles over $C$ of fixed type. 
Then there exists an irreducible family of
extensions of vector bundles, 
\begin{equation*}
\{ \ses{\mathcal{G}'_{z}}{\mathcal{F}'_{z}}{\mathcal{E}'_z}:z\in Z \}
\end{equation*}
such that every vector bundle $\mathcal{G}'_{z}$ is
isomorphic to some $\mathcal{G}_{x}$
and every extension
$\ses{\mathcal{G}_{x}}{\mathcal{F}}{\mathcal{E}_{y}}$
is isomorphic to one in this family.
\end{lemma}

\begin{proof}
After twisting by a suitable line bundle
of positive degree, we may assume that
$\Ext^{1}(\mathcal{O},\mathcal{G}_{x})=0$, for all $x\in X$,
and that every $\mathcal{E}_{y}$ is generated by global sections.
Suppose that each $\mathcal{E}_{y}$ is of type $(n,d)$.
Extending the usual dimension counting argument
in the Grassmannian $\Gr{n}{H^0(\mathcal{E}_{y})}$,
we may choose $n$ sections of $\mathcal{E}_{y}$ so that the induced
map $\rho:\mathcal{O}^n\to \mathcal{E}_{y}$ is an isomorphism of the fibres
at the general point of $C$ and drops rank by only 1 at other points.
Thus the cokernel of $\rho$ is the structure sheaf $\mathcal{T}_\xi$
of a subscheme $\xi$ of degree $d$ in $C$, that is,
$\mathcal{E}_{y}$ is an extension of $\mathcal{T}_\xi$ on top of
$\mathcal{O}^n$.

The parameter space of such subschemes $\xi$ is the $d$-fold symmetric 
product $C^{(d)}$, which is an irreducible
algebraic variety and which carries a universal family $\mathcal{T}$.
Since $\mathcal{T}_\xi$ is torsion, 
$\Hom(\mathcal{T}_{\xi},\mathcal{G}_{x}\oplus\mathcal{O}^{n})=0$
for all $\xi\in C^{(d)}$ and all $x\in X$.
Hence there is a vector bundle $\lambda:Z\rightarrow X\times C^{(d)}$ 
whose fibre above the point $(x,\xi)$ is
$\Ext(\mathcal{T}_{\xi},\mathcal{G}_{x}\oplus\mathcal{O}^{n})$ and this
carries a tautological family of extensions
\begin{equation*}
\{\ses{\mathcal{G}_{\pi_{1}\lambda(z)}\oplus\mathcal{O}^{n}}
{\mathcal{F}'_{z}}{\mathcal{T}_{\pi_{2}\lambda(z)}}:z\in Z \}.
\end{equation*}
Letting $\mathcal{G}'_z=\mathcal{G}_{\pi_{1}\lambda(z)}$ and
$\mathcal{E}'_z=\mathcal{F}'_{z}/\mathcal{G}'_{z}$,
we may replace $Z$ by the non-empty open set on which 
$\mathcal{F}'_{z}$ and $\mathcal{E}'_{z}$ are vector bundles
and obtain the required irreducible familty of extensions of vector bundles.
To see that every possible extension of $\mathcal{E}_{y}$
on top of $\mathcal{G}_{x}$ occurs note that every such extension 
has a 3 step filtration with $T_{\xi}$ on top of $\mathcal{O}^{n}$
on top of $\mathcal{G}_{x}$.
But, since $\Ext^{1}(\mathcal{O},\mathcal{G}_{x})=0$, the
extension at the bottom of this filtration splits and so it is
simply an extension of $T_{\xi}$ on top of 
$\mathcal{G}_{x}\oplus\mathcal{O}^{n}$.
\end{proof}

Using these lemmas, we have the following result.

\begin{proposition}\label{prop:EF}
For any type $\alpha=(r,d)$, let $h=\hcf{r}{d}$.
Then there is a unique type $\beta=(s,e)$ satisfying
\begin{enumerate}
 \item[(i)] $\eulchar{\beta}{\alpha}=h$,
 \item[(ii)] $r/h<s<2r/h$, if $h<r$, or $s=2$, if $h=r$. 
\end{enumerate}
Then, there exists a vector bundle $F$ of type $\beta$
such that for a general $E$ of type $\alpha$,
\begin{enumerate}
 \item[(iii)] $\hom(F,E)=h$ and $\ext(F,E)=0$,
 \item[(iv)]  the natural map 
   $\eval{F}{E}\colon\Hom(F,E)\otimes_{k}F\to E$ is surjective,
 \item[(v)] the bundle $E_1=\ker\eval{F}{E}$ is general
 and has $ext(E_1,F)=0$.
\end{enumerate}
\end{proposition}

\begin{proof}
To solve (i) we simply need to solve $sd-tr=h$ and set $e=t-(g-1)s$.
Given one solution $(s,t)$, the complete set of solutions is
$ \{ (s,t)+k(r/h,d/h) : k\in\ZZ \}$
which contains precisely one solution in the range (ii).
Part (iii) is provided by Lemma~\ref{lem:hirsch}.

For the main part of the proof,
the first step is to construct a short exact sequence
\begin{equation}
\label{ses:allh}
\ses{E_1}{F^h}{E}
\end{equation}
with $E$ of type $\alpha$, $F$ of type $\beta$ and
$\Ext(E_{1},F)=0=\Ext(F,E)$.

First suppose that $h=1$.
Then Lemma~\ref{lem:hirsch} implies that for generic $F$ and $E$
we have $\Ext(F,E)=0$ and, since $r<s$, the generic map is surjective.
Let $F'\to E'$ be a particular choice of such generic bundles and map
and let $E'_1$ be the kernel.
At this stage, we have $\Ext(F',E')=0$, but may not have $\Ext(E'_{1},F')=0$. 
On the other hand, 
\begin{equation*}
\eulchar{\beta-\alpha}{\beta}=\eulchar{\beta}{\alpha}+
\eulchar{\beta-\alpha}{\beta-\alpha}-\eulchar{\alpha}{\alpha}\geq
\eulchar{\beta}{\alpha}=1
\end{equation*}
since $s-r\leq r$.
Hence, Lemma~\ref{lem:hirsch} also implies that, for generic
$E_1$ and $F$ of types $\beta-\alpha$ and $\beta$ respectively,
$\Ext(E_1,F)=0$ and, since $s-r<s$, the generic map is an 
injection of vector bundles.
Let $E''_1\to F''$ be a particular choice of such generic bundles
and map and let $E''$ be its cokernel.
This time, we have $\Ext(E''_1,F'')=0$, but may not have
$\Ext(F'',E'')=0$.
 
But now we may include $E'_1$ and $E''_1$ in an 
irreducible family $\{ \mathcal{E}_{1,x} : x\in X \}$
by \cite{NarRam75} Proposition 2.6.
Then, by Lemma~\ref{lem:genext},
there is an irreducible family of extensions 
\begin{equation*}
\{ \ses{\mathcal{G}'_{z}}{\mathcal{F}'_{z}}{\mathcal{E}'_z}:z\in Z \}
\end{equation*}
which includes both $\ses{E'_1}{F'}{E'}$ and $\ses{E''_1}{F''}{E''}$.
Hence we may choose for (\ref{ses:allh}) a general extension in this family
and both $\Ext$ groups will vanish as required.

For an arbitrary value of $h$, we may obtain a sequence of the form
(\ref{ses:allh}) by taking the direct sum of $h$ copies
of one for $\alpha=\alpha/h$.

For the second step, suppose that we have a sequence of the form
(\ref{ses:allh}).
By \cite{NarRam75} Proposition 2.6, we may include 
$E_1$ in an irreducible family $\{\mathcal{E}_{1,x} : x\in X\}$,
whose generic member is general and for which every member satisfies
$\ext(\mathcal{E}_{1,x},F)=0$.  There is then a vector bundle 
$\lambda:Y\to X$
whose fibre at $x$ is $\Hom(\mathcal{E}_{1,x},F^h)$ and 
over $Y$ there is a 
tautological map $f_y:\mathcal{E}_{1,\lambda(y)}\to F^h$.
Replacing $Y$ by the non-empty open set on which $f_y$ is injective,
we have $\mathcal{E}_y=\coker f_y$ of type $\alpha$.
We may further replace $Y$ by the non-empty open set
on which $\ext(F,\mathcal{E}_{y})=0$.

Now observe that the homomorphism $F^{h}\to\mathcal{E}_{y}$ 
must be isomorphic to
$\eval{F}{\mathcal{E}_{y}}\colon\Hom(F,\mathcal{E}_{y})\otimes_{k}F\to
\mathcal{E}_{p}$, because a linear dependence between the $h$ components
of the homomorphism from $F^{h}$ to $\mathcal{E}_{y}$ would imply that
$F$ is a summand of the kernel, 
which would contradict $\ext(\mathcal{E}_{1,x},F)=0$.

If we consider just the family $\{\mathcal{E}_y:y\in Y\}$,
then, as before,
we may include this in an irreducible family 
$\{\mathcal{E}_z : z\in Z \}$,
whose generic member is general 
and such that $\ext(F,\mathcal{E}_{z})=0$ and 
$\hom(F,\mathcal{E}_{z})=h$ for every $z\in Z$.
Hence the kernel $\mathcal{E}'_{1,z}$ of the homomorphism
$\eval{F}{\mathcal{E}_{z}}\colon\Hom(F,\mathcal{E}_{z})\otimes_{k}F\to
\mathcal{E}_{z}$ is general and satisfies 
$\ext(\mathcal{E}'_{1,z},F)=0$,
because this was already true over $Y$.
Thus we have all the properties we require. 
\end{proof}

Proposition~\ref{prop:EF} shows that we have a dominant rational map
$$\lambda_F\colon \moduli{r}{d}\ramto\moduli{r_{1}}{d_{1}}:[E]\mapsto [E_1],$$ 
where $E_1=\ker\eval{F}{E}$ has type $(r_1,d_1)=h(s,e)-(r,d)$.
One may immediately check the following. 
\begin{lemma}
 \label{lem:r1d1}
  The type $(r_1,d_1)$ of $E_1$ satisfies
  \begin{enumerate}
   \item[(i)] if $h<r$, then $r_1<r$,
   \item[(ii)] $h_1=\hcf{r_1}{d_1}$ is divisible by $h$,
   \item[(iii)] $\det(E_1)\isom\det(F)^h\det(E)^{-1}$.
  \end{enumerate}
\end{lemma}

The proof of Proposition~\ref{prop:EF} shows
 that the fibre of $\lambda_F$ above a closed 
point $[E_1]$ is birationally the Grassmannian
of $h$-dimensional subspaces of $\Hom(E_1,F)$.
However, this bundle of Grassmannians may be `twisted',
that is, it may not be locally trivial in the Zariski topology.
In fact, it will fail to be locally trivial whenever $h\neq 1$ and
will not be birationally linear whenever $h_1\neq h$,
but we will be able to measure how twisted it is
using a Brauer class on $\moduli{r_1}{d_1}$ and then
compare $\lambda_F$ to another Grassmannian bundle with
the same twisting, but smaller fibres, to construct inductively
our birationally linear map.

In fact, Proposition~\ref{prop:EF} also provides us with the 
way of constructing this Brauer class, 
because it yields a description of $\moduli{r}{d}$ as a quotient
of an open set in the quot scheme $\quot{F^h}{r}{d}$ by $\PGL{h}$.
We describe this in detail in the next section. 

\section{Brauer classes and free $\PGL{}$ actions} \label{s2}

In this section, we collect a number of results about free actions
of the projective general linear group $\PGL{}$,
which allow us to define and
compare the Brauer classes we are interested in.

Recall that the Brauer group of a field $k$ may be described as
consisting of classes represented by central simple algebras $A$ over
the field and that $[A_1]=[A_2]$ in the Brauer group if and only if
$A_1$ and $A_2$ are Morita equivalent or equivalently
$A_{1}^{o}\otimes A_{2}$ is isomorphic to $M_{n}(k)$ for a suitable
integer $n$ where $A^{o}$ is the opposite algebra to $A$.  This is
equivalent to saying that there is an $A_{1}$, $A_{2}$ bimodule of
dimension $n$ where $n^{2}=\dim A_{1}\dim A_{2}$. The product in the
Brauer group is induced by the tensor product of algebras.

The reader may wish to consult \cite{Draxl} for further discussion of
the Brauer group and central simple algebras.

\begin{definition}
Let $X$ be an affine algebraic variety on which the algebraic group $\PGL{n}$
acts freely.
Over the quotient variety $X/\PGL{n}$ there is a bundle of
central simple algebras $M_{n}(k)\times^{\PGL{n}}X$
of dimension $n^{2}$.
At the generic point, this is a central simple algebra over
the function field $\funfield{X/\PGL{n}}$
and hence defines a class in the Brauer group of $\funfield{X/\PGL{n}}$.
We shall denote this class by $\Br{X/\PGL{n}}{X}{\sigma}$.
\end{definition}

It is important to note that $\Br{X/\PGL{n}}{X}{\sigma}$
depends on the action of $\PGL{n}$ on $X$ and 
not just on the quotient space $Y=X/\PGL{n}$.
Note also that the bundle of central simple algebras
$B=M_{n}(k)\times^{\PGL{n}}X$ over $Y$ 
is essentially equivalent to the $\PGL{n}$ action on $X$,
because $X$ can be recovered as the $Y$-scheme that represents the
functor of isomorphisms between $B$ and the trivial bundle of central
simple algebras $M_{n}(k)\times Y$ over $Y$.
The $\PGL{n}$ action is recovered via its action on $M_{n}(k)$.
Moreover, we have the following.
\begin{lemma}
\label{t0}
Let $\PGL{n}$ act freely on affine algebraic varieties $X_{1}$ and $X_{2}$.
Let
$$\phi\colon X_{1}/\PGL{n}\rightarrow X_{2}/\PGL{n}$$
be a dominant rational map.
Then there is a $\PGL{n}$-equivariant
dominant rational map $\Phi\colon X_{1}\rightarrow X_{2}$
making the following diagram commute
\begin{equation*}
\begin{CD}
X_{1}@>\Phi>>X_{2}\\
@VVV         @VVV\\
X_{1}/\PGL{n}@>\phi>>X_{2}/\PGL{n} 
\end{CD}
\end{equation*} if and only if
$\Br{X_{1}/\PGL{n}}{X_{1}}{\sigma_{1}}
=\phi^{-1}\Br{X_{2}/\PGL{n}}{X_{2}}{\sigma_{2}}$.
\end{lemma}
\begin{proof}
After restricting to suitable open subvarieties and taking the
pullback along $\phi$ we may assume that $\phi$ is the identity
map. We have two distinct $\PGL{n}$ bundles. These have the same Brauer
class if and only if over a suitable open subvariety of $X/\PGL{n}$ the
associated bundles of central simple algebras are isomorphic or
equivalently the two $\PGL{n}$ bundles are isomorphic over this open
subvariety.
\end{proof}

We can now define the Brauer classes on (the function fields of)
our moduli spaces that we will use in the rest of the paper.
For each type $(r,d)$, fix one vector bundle $F$,
which is general in the sense of Proposition~\ref{prop:EF} and
recall that $h=\hcf{r}{d}$. 
Let $X_{r,d}$ be the open subset of $\quot{F^h}{r}{d}$,
which parametrizes (up to scaling) quotients $q\colon F^h\to E$ 
of type $(r,d)$ which are stable bundles and for which the induced 
map $k^{h}\to\Hom(F,E)$ is an isomorphism. 
The obvious action of $GL_{h}=\Aut(F^h)$ induces a free action of
$\PGL{h}$ on $X_{r,d}$ and the map $X_{r,d}\to \moduli{r}{d}$, 
which forgets the quotient map $q$, identifies
$X_{r,d}/\PGL{h}$ with an open dense subset of $\moduli{r}{d}$
and, in particular, identifies their function fields. Since
$\moduli{r}{d}$ is a projective variety we may replace $X_{r,d}$ by an
open dense affine $\PGL{h}$-equivariant subset of itself by taking the
inverse image of some open dense affine subset of $\moduli{r}{d}$
contained in the image of $X_{r,d}$.

\begin{definition}\label{d2}
For every type $(r,d)$, the Brauer class $\brauer{r}{d}$
on $\moduli{r}{d}$ is the class corresponding to
$\Brr{X_{r,d}/\PGL{h}}$ after we identify
$\funfield{X_{r,d}/\PGL{h}}$ with $\funfield{\moduli{r}{d}}$
as described above.
\end{definition}

There are more general Brauer classes that arise naturally on $X/\PGL{n}$,
which we now describe and relate to $\Brr{X/\PGL{n}}$.
Let $P$ be a vector bundle over the algebraic variety $X$ on which
$\GL{n}$ acts lifting the action of $\PGL{n}$ on $X$ such that
$k^{\ast}$ acts with weight $w$ on the fibres of $P$.
We will call such a bundle $P$ a vector bundle
of weight $w$ on $X$; the $\GL{n}$ action on $P$ lifting
the $\PGL{n}$ action on $X$ will be implicit.
If $P$ is a vector bundle of weight $0$, then
$\PGL{n}$ acts on $P$ and $P/\PGL{n}$ is a vector bundle over
$X/\PGL{n}$. If $P$ is a vector bundle of weight $w$ then 
$P\dual\otimes P$ is a vector bundle of weight $0$ and 
$P\dual\otimes P/\PGL{n}$ is a bundle of central simple algebras over
$X/\PGL{n}$. The bundle of central simple algebras associated to the
$\PGL{n}$ action of $X$ is the special case where $P$ is taken to be
the bundle of weight $1$ over $X$ given by $k^{n}\times X$,
with $\GL{n}$ acting diagonally,
since we may identify $M_n(k)$ with $(k^{n})\dual\otimes k^{n}$.
We define the Brauer class defined by $P$ to be the Brauer
class of the central simple algebra over $\funfield{X/\PGL{n}}$
defined by the generic fibre of the bundle of central simple algebras 
$P\dual\otimes P/\PGL{n}$.

\begin{lemma}
\label{t1}
Let $P$ be a vector bundle of weight $w$ over an algebraic variety
$X$ on which $\PGL{n}$ acts freely.
Then the Brauer class defined by $P$ is
$w\Br{X/\PGL{n}}{X}{\sigma}$.
\end{lemma}
\begin{proof}
Let $P$ and $Q$ be vector bundles of weight $w$.
Then $P\dual\otimes Q$ is a vector bundle of weight $0$ 
and $P\dual\otimes Q/\PGL{n}$ is a
vector bundle over $X/\PGL{n}$. It has a structure of a bimodule
with $P\dual\otimes P/\PGL{n}$ acting on the left and
$Q\dual\otimes Q/\PGL{n}$ acting on the right.
Over the generic point of $X/\PGL{n}$ it defines a Morita equivalence
between (the generic fibres of) $P\dual\otimes P/\PGL{n}$ and 
$Q\dual\otimes Q/\PGL{n}$.
Hence the Brauer classes defined by $P$ and $Q$ are equal;
in other words the Brauer class depends only on the weight. 

Now, if $w>0$, then $Q_{w}=\left(k^{n}\right)^{\otimes w}\times X$ with the
diagonal action of $\GL{n}$ is a vector bundle of weight $w$ and
$Q_{w}\dual\otimes Q_{w}/\PGL{n}$ is the $w$th tensor power of
$Q_{1}\dual\otimes Q_{1}/\PGL{n}$.
Since the class defined by $Q_{1}$ is $\Br{X/\PGL{n}}{X}{\sigma}$,
the class defined by $Q_w$ is $w\Br{X/\PGL{n}}{X}{\sigma}$.  
On the other hand,
$Q_{-w}=\left((k^{n})\dual\right)^{\otimes w}\times X$
 with the diagonal action of
$\GL{n}$ is a vector bundle of weight $-w$.
In particular, $Q_{-1}\dual\otimes Q_{-1}/\PGL{n}$
is the sheaf of algebras opposite to $Q_{1}\dual\otimes Q_{1}/\PGL{n}$
and therefore the class defined by $Q_{-1}$ is
$-\Br{X/\PGL{n}}{X}{\sigma}$ and, as above, the class defined by
$Q_{-w}$ is $-w\Br{X/\PGL{n}}{X}{\sigma}$.
Finally, $\cO_X$ is a vector bundle of weight $0$
and defines the class 0.
\end{proof}

Thus, if $P$ is a vector bundle of weight $1$ and rank $r$,
then the Brauer class $\Br{X/\PGL{n}}{X}{\sigma}$ is represented
by a central simple algebra of dimension $r^2$, 
namely $P\dual\otimes P/\PGL{n}$.
It will be important to observe that, birationally, the converse is true.
More precisely, we have the following.

\begin{lemma}
Let $\PGL{n}$ act freely on an algebraic variety $X$ and suppose that
the Brauer class $w\Br{X/\PGL{n}}{X}{\sigma}$ is represented by a
central simple algebra $S$ of dimension $s^2$ over $\funfield{X/\PGL{n}}$.
Then there exists a $\PGL{n}$-equivariant open subset
$Y$ of $X$ and a vector bundle $Q$ of weight $w$ over $Y$ whose rank
is $s$.
\end{lemma}
\begin{proof}
Let $P$ be a vector bundle of weight $w$ and rank $p$.
It is enough to deal with the case where $S$ is a division
algebra since the remaining cases are all matrices over this and hence
the values for $s$ that arise are all multiples of this. In
particular, therefore, we may assume that $s$ divides $p$. If $s=p$,
there is nothing to prove so we may assume that $s<p$. Thus at the
generic point of $P\dual\otimes P/\PGL{n}$, there is an idempotent of
rank $s$.  This idempotent is defined over some open subset of $X$
which is $\PGL{n}$-equivariant since the idempotent is $\PGL{n}$-invariant
and gives a decomposition $P\cong P_1\oplus P_2$ as a direct sum of
vector bundles which are $\GL{n}$-equivariant one of which has rank
$s$. These bundles have weight $w$ since they are subbundles of $P$
which has weight $w$.
\end{proof}

We now come to the main object of this section, to describe the
relationship between the Brauer classes considered above and twisted
Grassmannian bundles such as
$\lambda_F\colon\moduli{r}{d}\to\moduli{r_{1}}{d_{1}}$. 
We start in
the general context of Grassmannian bundles associated to a vector
bundle $P$ of weight $w$, although in the end we will only need to
consider weights $\pm 1$.  Let $j<\rk(P)$ be a positive integer.  Then
$\PGL{n}$ acts freely on the bundle of Grassmannians $\Gr{j}{P}$ over
$X$ and $\phi\colon \Gr{j}{P}/\PGL{n}\rightarrow X/\PGL{n}$ is a
Grassmannian bundle over $X/\PGL{n}$ that is usually not trivial in
the Zariski topology.  Since the map from $\Gr{j}{P}$ to $X$ is
$\PGL{n}$-equivariant the Brauer class
$\Br{\Gr{j}{P}/\PGL{n}}{\Gr{j}{P}}{\sigma}$ is just the pullback of
the Brauer class $\Br{X/\PGL{n}}{X}{\sigma}$. We can also realise the
algebraic variety $\Gr{j}{P}/\PGL{n}$ as a quotient variety for a free
action of the algebraic group $\PGL{j}$ on the partial frame bundle of
$j$ linearly independent sections of the vector bundle $P$ and we will
need to know how to relate the two Brauer classes we obtain in this
way.

We must take care to differentiate two distinct ways of
constructing the partial frame bundle.
Let $S$ be the universal sub-bundle on $\Gr{j}{P}$.
Let $\Fr{j}{P}$ be the `covariant' partial frame bundle,
whose fibre at $x$ consists of isomorphisms $k^{j}\to S_x$
and let $\Frdual{j}{P}$ be the `contravariant' partial frame bundle,
whose fibre at $x$ consists of isomorphisms $(k^{j})\dual\to S_x$.
Then $\GL{j}$ acts freely on both $\Fr{j}{P}$ and $\Frdual{j}{P}$
and the quotient variety is $\Gr{j}{P}$ in both cases.
The difference is that the pullback of $S$ to $\Fr{j}{P}$
is the trivial bundle with fibre $k^{j}$ on which $\GL{j}$
acts with weight $1$, while the pullback of $S$ to $\Frdual{j}{P}$
is the trivial bundle with fibre $(k^{j})\dual$ on which $\GL{j}$
acts with weight $-1$.
The obvious isomorphism between the two frame bundles is compatible
with the transpose inverse automorphism of $\GL{j}$, but not
with the identity automorphism.

The action of $\GL{n}$ lifts from $\Gr{j}{P}$ to $\Fr{j}{P}$
and $\Frdual{j}{P}$, so both carry an action
of $\GL{j}\times \GL{n}$.
The kernel of each action is isomorphic to $k^\ast$,
but in the covariant case it is $\{(t^{w}I,t I):t\in k^\ast\}$,
while in the contravariant case it is 
$\{(t^{w}I,t^{-1} I):t\in k^*\}$.
(Recall that $w$ is the weight of the action of $\GL{n}$ on $P$.)
Hence, both $\Fr{j}{P}/\GL{n}$ and $\Frdual{j}{P}/\GL{n}$
carry free actions of $\PGL{j}$ which determine Brauer classes
on the quotient, which is equal to $\Gr{j}{P}/\PGL{n}$ in both cases.

We summarise the maps considered above in the following commutative diagram
for the case of the covariant partial frame bundle.
Note that the groups that appear as labels below the arrows indicate that the 
maps are quotient maps by a (generically) free action of the group.
\begin{equation}\label{diag:GrFr}
 \setlength{\unitlength}{55pt}
  \begin{picture}(3,3)(0,0)
   \put(1,0){\object{$X/\PGL{n}$}}
   \put(0,1){\object{$X$}}
   \put(0.25,0.75){\vector(1,-1){0.5}}
   \put(0.4,0.5){\nwlabel{$\PGL{n}$}}
   \put(2,1){\object{$\Gr{j}{P}/PGL{n}$}}
   \put(1.75,0.75){\vector(-1,-1){0.5}}
   \put(1.5,0.5){\nelabel{$\phi$}}
   \put(1,2){\object{$\Gr{j}{P}$}}
   \put(0.75,1.75){\vector(-1,-1){0.5}}
   \put(1.25,1.75){\vector(1,-1){0.5}}
   \put(1.4,1.5){\nwlabel{$\PGL{n}$}}
   \put(3,2){\object{$\Fr{j}{P}/\GL{n}$}}
   \put(2.75,1.75){\vector(-1,-1){0.5}}
   \put(2.6,1.5){\swlabel{$\PGL{j}$}}
   \put(2,3){\object{$\Fr{j}{P}$}}
   \put(1.75,2.75){\vector(-1,-1){0.5}}
   \put(1.5,2.5){\swlabel{$\GL{j}$}}
   \put(2.25,2.75){\vector(1,-1){0.5}}
   \put(2.5,2.5){\nwlabel{$\GL{n}$}}
  \end{picture}
\end{equation}
The diagram for the contravariant frame bundle is of the identical form,
but the need to distinguish the two cases is made clear by the 
following result, which describes the relationship between the
Brauer classes determined by all the $\PGL{}$ actions in the diagram.

\begin{lemma}
\label{t2}
Let $P$ be a bundle of weight $w$ on $X$.
Then, in the notation described above,
\begin{align*}
\Brr{\Gr{j}{P}/\PGL{n}}&=\phi^*\left(\Brr{X/\PGL{n}}\right)\\
\Brr{(\Fr{j}{P}/\GL{n})\big/\PGL{j}} &= w\Brr{\Gr{j}{P}/\PGL{n}}\\
\Brr{(\Frdual{j}{P}/\GL{n})\big/\PGL{j}} &= -w\Brr{\Gr{j}{P}/\PGL{n}}
\end{align*}
\end{lemma}

\begin{proof}
The first equality follows immediately from the fact that the lower
diamond in (\ref{diag:GrFr}) is an equivariant pullback.
The action of $\GL{n}$ on $P$ over $X$ lifts naturally to an action of
$\GL{n}$ on the universal subbundle $S$ over $\Gr{j}{P}$.
Hence $S$ has weight $w$ and so the Brauer class
on $\Gr{j}{P}/\PGL{n}$ represented by $S\dual\otimes S/\PGL{n}$
is $w\Brr{\Gr{j}{P}/\PGL{n}}$ by Lemma~\ref{t1}.
As already observed, the pullback $S'$ of $S$ to $\Fr{j}{P}$
is trivial with fibre $k^j$ and so the quotient by $\GL{n}$
also gives a trivial bundle $S''$ with fibre $k^j$ on $\Fr{j}{P}/\GL{n}$.
But then $\Brr{(\Fr{j}{P}/\GL{n})\big/\PGL{j}}$
is equal to the Brauer class represented by 
$(S'')\dual\otimes(S'')/\PGL{j}$, which is equal to
the Brauer class represented by $S\dual\otimes S/\PGL{n}$,
completing the proof in the covariant case.
The proof in the contravariant case is identical except that
now $S'$ and $S''$ are trivial with fibre $(k^j)\dual$ so that
$\Brr{(\Frdual{j}{P}/\GL{n})\big/\PGL{j}}$ is the negative of
the class represented by $(S'')\dual\otimes(S'')/\PGL{j}$.
\end{proof}

We may now describe the rational map 
$\lambda_F\colon\moduli{r}{d}\to\moduli{r_{1}}{d_{1}}$
as a Grassmannian bundle of the type described above
and determine the behaviour of the Brauer classes under pullback.
Recall that there exists a vector bundle $F_{1}$
and an open subset $X_{r_{1},d_{1}}$ of
$\quot{F_{1}^{h_{1}}}{r_{1}}{d_{1}}$ such that the map to
$\moduli{r_{1}}{d_{1}}$ is birational to the 
$\PGL{h_{1}}$ quotient map. 

\begin{proposition}
\label{t10}
On an open subset of $X_{r_{1},d_{1}}$, there exists a vector
bundle $P$ of weight $-1$
and of rank $lh_{1}$ for some integer $l$
such that the rational map 
$$\lambda_F\colon\moduli{r}{d}\to\moduli{r_{1}}{d_{1}}$$
is birational to the Grassmannian bundle 
$$\phi:\Gr{h}{P}/\PGL{h_{1}}\to X_{r_1,d_1}/\PGL{h_{1}}.$$
Furthermore,
\begin{equation*}
 \lambda_F^*(\brauer{r_{1}}{d_{1}})=\brauer{r}{d}.
\end{equation*}
\end{proposition}
\begin{proof}
The idea of the proof is that since $F$ is general, the set of
quotients $p:F \otimes V \to E \in X_{r,d}$ such that $E_1 :=\ker
p \in \moduli{r_1}{d_1}$, $Ext(F,E)=0$ and $Ext(E_1,F)=0$ is
not empty. It is bijective to the set of $h$-dimensional
subspaces $V^{\vee} \subset Hom(E_1,F)$ such that $p:E_1 \to F
\otimes V$ is injective, $E:=\coker p \in\moduli{r}{d}$, 
$Ext(E_1,F)=0$ and $Ext(F,E)=0$.
We fill in the details below. 

Consider the open set in $X_{r_{1},d_{1}}$ parametrizing those 
$q_1\colon F_{1}^{h_{1}}\to E_{1}$
for which $\Ext(E_{1},F)=0$.
Over this open set there is a vector bundle $P$ whose fibre at $[q_1]$
is $\Hom(E_{1},F)$.
Since $\PGL{h_1}$ acts with weight 1 on $E_1$,
it acts with weight $-1$ on $P$.
We claim that $\Gr{h}{P}/\PGL{h_{1}}$ is birational to $\moduli{r}{d}$.
To see this, consider the contravariant partial frame bundle 
$\Frdual{h}{P}$ whose fibre over $[q_1]\in X_{r_{1},d_{1}}$
is naturally identified with the set of maps $p\colon E_1\to F^h$
for which the induced map $(k^h)\dual\to\Hom(E_1,F)$ is injective.
On an open subset in $\Frdual{h}{P}$,
the map $p$ is injective as a map of bundles
and its cokernel $q:F^h\to E$ gives a point in $X_{r,d}$.
Since $p$, but not $q_1$, is determined by $q$
we see that $\Frdual{h}{P}/\GL{h_{1}}$ is birational
to $X_{r,d}$ and so $(\Frdual{h}{P}/\GL{h_{1}})/\PGL{h}$
is birational to $\moduli{r}{d}$.
Since $\Gr{h}{P}$ is $\Fr{h}{P}/\GL{h}$,
we deduce that $\Gr{h}{P}/\PGL{h_{1}}$ is birational to $\moduli{r}{d}$.

We can arrange all the rational maps we have considered above into
the following diagram of the form of (\ref{diag:GrFr}).
\begin{equation}\label{diag:lambdaF}
 \setlength{\unitlength}{55pt}
  \begin{picture}(3,3)(0,0)
   \put(1,0){\object{$\moduli{r_1}{d_1}$}}
   \put(0,1){\object{$X_{r_1,d_1}$}}
   \put(0.25,0.75){\vector(1,-1){0.5}}
   \put(0.4,0.5){\nwlabel{$\PGL{h_1}$}}
   \put(2,1){\object{$\moduli{r}{d}$}}
   \put(1.75,0.75){\vector(-1,-1){0.5}}
   \put(1.5,0.5){\nelabel{$\lambda_F$}}
   \put(1,2){\object{$\Gr{h}{P}$}}
   \put(0.75,1.75){\vector(-1,-1){0.5}}
   \put(1.25,1.75){\vector(1,-1){0.5}}
   \put(1.4,1.5){\nwlabel{$\PGL{h_1}$}}
   \put(3,2){\object{$X_{r,d}$}}
   \put(2.75,1.75){\vector(-1,-1){0.5}}
   \put(2.6,1.5){\swlabel{$\PGL{h}$}}
   \put(2,3){\object{$\Frdual{h}{P}$}}
   \put(1.75,2.75){\vector(-1,-1){0.5}}
   \put(1.5,2.5){\swlabel{$\GL{h}$}}
   \put(2.25,2.75){\vector(1,-1){0.5}}
   \put(2.5,2.5){\nwlabel{$\GL{h_1}$}}
  \end{picture}
\end{equation}
Since $P$ has weight $-1$, the first and last formulae in
Lemma~\ref{t2} give
\begin{equation*}
 \lambda_F^*(\brauer{r_1}{d_1})
 =\Brr{\Gr{h}{P}/\PGL{h_{1}}}
 =\brauer{r}{d}
\end{equation*}
which completes the proof.
\end{proof}

\section{The Hecke correspondence}\label{sec:Hecke}

One of the main ideas of the paper is to compare the 
(birationally) twisted Grassmannian bundle 
$\lambda_F\colon\moduli{r}{d}\to\moduli{r_{1}}{d_{1}}$
to another Grassmannian bundle which is twisted by the same amount
but has smaller fibres.
This second bundle is provided by the Hecke correspondence,
which we describe in this section.
Within this section, we may let $h$ and $h_1$ be arbitrary
integers with $h\leq h_1$. Only later, will we
need to use the fact that $h$ actually divides $h_1$.

Let $\parmod{h_{1}}{0}{h}$ be the moduli space of parabolic
bundles, which parametrizes pairs
consisting of a bundle (or locally free sheaf) $\mathcal{E}_1$ 
of type $(h_{1},0)$ together with a locally free subsheaf
$\mathcal{E}_2\subset\mathcal{E}_1$ such that the quotient 
$\mathcal{E}_1/\mathcal{E}_2$
is isomorphic to $\bigl(\mathcal{O}_x\bigr)^h$ for a fixed
point $x\in C$.
In order to specify a projective moduli space exactly,
we would need to specify parabolic weights to determine notions
of stability and semistability.
However, we are only interested in this space up to birational
equivalence and it is known (\cite{BodYok} Section 4)
that the birational type of the moduli space does not depend on
the choice of parabolic weights.
Indeed, we may choose to let $\parmod{h_{1}}{0}{h}$
denote the dense open set of quasi-parabolic bundles
$\mathcal{E}_2\subset\mathcal{E}_1$ that are stable for all choices
of parabolic weights. 

The type of $\mathcal{E}_2$ must be $(h_{1},-h)$ and there
are two dominant rational maps 
\begin{align*}
 \theta_1&\colon\parmod{h_{1}}{0}{h}\ramto\moduli{h_{1}}{0}
 \colon [\mathcal{E}_2\subset\mathcal{E}_1]\mapsto [\mathcal{E}_1]\\
 \theta_2&\colon\parmod{h_{1}}{0}{h}\ramto\moduli{h_{1}}{-h}
 \colon [\mathcal{E}_2\subset\mathcal{E}_1]\mapsto [\mathcal{E}_2]
\end{align*}

The key point is that, like $\lambda_F$, the
maps $\theta_1$ and $\theta_2$ are (birational to)
twisted Grassmannian bundles whose twisting is measured
by the Brauer classes $\brauer{h_{1}}{0}$ and $\brauer{h_{1}}{-h}$
respectively.
Furthermore, as we shall show below, these two Brauer classes
pull back to the same class on $\parmod{h_{1}}{0}{h}$.

To construct $\parmod{h_{1}}{0}{h}$ birationally from $\moduli{h_{1}}{0}$,
let $H_1$ be the vector bundle over $X_{h_1,0}$ whose fibre over 
the point $\left[q_1:F_1^{h_1}\to E_1\right]$ is $\Hom(E_1,\cO_x)$,
where $\cO_x$ is the structure sheaf of the point $x\in C$.
Then $H_1$ is a vector bundle of weight $-1$
and $\parmod{h_{1}}{0}{h}$ is birational to $\Gr{h}{H_1}/\PGL{h_1}$.
To see this, consider the contravariant frame bundle $\Frdual{h}{H_1}$.
A point in the fibre over $[q_1]$ may be identified with a map
\begin{equation}\label{eq:p1}
 p\colon E_1\rightarrow (\cO_x)^h
\end{equation}
such that the induced map $(k^h)\dual\to\Hom(E,\cO_x)$ is injective.
If we restrict to the open set on which $p$ is also surjective
so that it determines a quasi-parabolic structure,
then the map to $\parmod{h_{1}}{0}{h}$ which forgets $p$
and $q_1$ is precisely the quotient by $\GL{h}$ that gives
$\Gr{h}{H_1}$, followed by the quotient by $\PGL{h_1}$.

Thus we have another diagram of the form of (\ref{diag:GrFr}).
\begin{equation}\label{diag:theta1}
 \setlength{\unitlength}{55pt}
  \begin{picture}(3,3)(0,0)
   \put(1,0){\object{$\moduli{h_1}{0}$}}
   \put(0,1){\object{$X_{h_1,0}$}}
   \put(0.25,0.75){\vector(1,-1){0.5}}
   \put(0.4,0.5){\nwlabel{$\PGL{h_1}$}}
   \put(2,1){\object{$\parmod{h_{1}}{0}{h}$}}
   \put(1.75,0.75){\vector(-1,-1){0.5}}
   \put(1.5,0.5){\nelabel{$\theta_1$}}
   \put(1,2){\object{$\Gr{h}{H_1}$}}
   \put(0.75,1.75){\vector(-1,-1){0.5}}
   \put(1.25,1.75){\vector(1,-1){0.5}}
   \put(1.4,1.5){\nwlabel{$\PGL{h_1}$}}
   \put(3,2){\object{$\Frdual{h}{H_1}/\GL{h_1}$}}
   \put(2.75,1.75){\vector(-1,-1){0.5}}
   \put(2.6,1.5){\swlabel{$\PGL{h}$}}
   \put(2,3){\object{$\Frdual{h}{H_1}$}}
   \put(1.75,2.75){\vector(-1,-1){0.5}}
   \put(1.5,2.5){\swlabel{$\GL{h}$}}
   \put(2.25,2.75){\vector(1,-1){0.5}}
   \put(2.5,2.5){\nwlabel{$\GL{h_1}$}}
  \end{picture}
\end{equation}
Hence, by Lemma~\ref{t2}, we have
\begin{equation}\label{eq:theta1}
\theta_1^*(\brauer{h_1}{0})=\Brr{(\Frdual{h}{H_1}/\GL{h_1})/\PGL{h}}
\end{equation}

Now we construct $\parmod{h_{1}}{0}{h}$
birationally from $\moduli{h_{1}}{-h}$.
To preserve the generality of this section, let $m=\hcf{h_1}{h}$,
but note that $m=h$ in the case of real interest.
Let $H_2$ be the vector bundle over $X_{h_1,-h}$
whose fibre above a point $\left[q_2:F_2^m\to E_2\right]$ is 
$\Ext(\cO_x,E_2)$.
Then $H_2$ has weight $1$ and $\parmod{h_{1}}{0}{h}$
is also birational to $\Gr{h}{H_2}/\PGL{m}$.
This follows, as above, by considering the open set
in $\Fr{h}{H_2}$ parametrizing extensions 
\begin{equation}\label{eq:p2}
\ses{E_2}{E_1}{(\cO_x)^h}.
\end{equation}
such that the induced map $k^h\to \Ext(\cO_x,E_2)$ is injective.
The moduli space $\parmod{h_{1}}{0}{h}$ arises (birationally)
by taking the quotient by $\GL{h}$ and then $\PGL{m}$.

Thus, again, we have a diagram of the form of (\ref{diag:GrFr}).
\begin{equation}\label{diag:theta2}
 \setlength{\unitlength}{55pt}
  \begin{picture}(3,3)(0,0)
   \put(1,0){\object{$\moduli{h_1}{-h}$}}
   \put(0,1){\object{$X_{h_1,-h}$}}
   \put(0.25,0.75){\vector(1,-1){0.5}}
   \put(0.4,0.5){\nwlabel{$\PGL{m}$}}
   \put(2,1){\object{$\parmod{h_{1}}{0}{h}$}}
   \put(1.75,0.75){\vector(-1,-1){0.5}}
   \put(1.5,0.5){\nelabel{$\theta_2$}}
   \put(1,2){\object{$\Gr{h}{H_2}$}}
   \put(0.75,1.75){\vector(-1,-1){0.5}}
   \put(1.25,1.75){\vector(1,-1){0.5}}
   \put(1.4,1.5){\nwlabel{$\PGL{m}$}}
   \put(3,2){\object{$\Fr{h}{H_2}/\GL{m}$}}
   \put(2.75,1.75){\vector(-1,-1){0.5}}
   \put(2.6,1.5){\swlabel{$\PGL{h}$}}
   \put(2,3){\object{$\Fr{h}{H_2}$}}
   \put(1.75,2.75){\vector(-1,-1){0.5}}
   \put(1.5,2.5){\swlabel{$\GL{h}$}}
   \put(2.25,2.75){\vector(1,-1){0.5}}
   \put(2.5,2.5){\nwlabel{$\GL{m}$}}
  \end{picture}
\end{equation}
and, by Lemma~\ref{t2}, we have
\begin{equation}\label{eq:theta2}
\theta_2^*(\brauer{h_1}{-h})=\Brr{(\Fr{h}{H_2}/\GL{m})/\PGL{h}}
\end{equation}

But now we simply need to observe that the data in (\ref{eq:p1})
and in (\ref{eq:p2}) have the same form and differ only in the imposition
of different open conditions.
Thus we may identify open subsets of $\Fr{h}{H_2}/\GL{m}$
and $\Fr{h}{H_1}/\GL{h_1}$.
One could in principle identify both of these as open subsets
of an appropriate fine moduli space for such data.
Furthermore, this identification is compatible with the
$\PGL{h}$ actions and so we may combine (\ref{eq:theta1}) and
(\ref{eq:theta2}) to obtain
\begin{equation}\label{eq:thetas}
  \theta_2^*(\brauer{h_1}{-h})=\theta_1^*(\brauer{h_1}{0}).
\end{equation}

\section{Construction of birationally linear maps} \label{s11}

We may now proceed with the proof of the main theorem of the paper
on the existence of a birationally linear map from 
the moduli space $\moduli{r}{d}$ of vector bundles of 
rank $r$ and degree $d$ to
the moduli space $\moduli{h}{0}$.
The proof goes by induction on the stronger statement that 
there is such a birationally linear rational map
that preserves the Brauer classes defined 
in the Section~\ref{s2}.
More precisely, we prove the following theorem.

\begin{theorem}
\label{thm:main}
Let $\brauer{r}{d}$ be the Brauer class on $\moduli{r}{d}$
defined for every type $(r,d)$ in Definition~\ref{d2}
and let $h=\hcf{r}{d}$.
Then there exists a birationally linear map 
$\mu\colon \moduli{r}{d}\ramto\moduli{h}{0}$
such that $\mu^*(\brauer{h}{0})=\brauer{r}{d}$.
\end{theorem}
\begin{proof}
If $r$ divides $d$, then $h=r$ and $\moduli{r}{d}$ is isomorphic 
to $\moduli{h}{0}$ by tensoring with a line bundle of degree $d/r$.
This isomorphism may be taken to be $\mu$ and it preserves the Brauer class.
Otherwise, we saw in Section~\ref{sec:step1} how to construct a map
$\lambda_F\colon\moduli{r}{d}\ramto\moduli{r_1}{d_1}$ with $r_1<r$
and we proved in Proposition~\ref{t10} that 
\begin{equation}
 \lambda_F^*(\brauer{r_{1}}{d_{1}})=\brauer{r}{d}.
\end{equation}
We construct the map $\mu$, by induction on the rank $r$,
as the composite of the top row of the following commutative
diagram of dominant rational maps which combines $\lambda_F$
and the Hecke correspondence described in Section~\ref{sec:Hecke}.
The other elements of the diagram we will explain next.
\begin{equation}\label{eq:mainpic}
 \setlength{\unitlength}{80pt}
  \begin{picture}(4,1.2)(0,-0.1)
   \put(0,1){\object{$\moduli{r}{d}$}}
   \put(0.25,0.75){\vector(1,-1){0.5}}
   \put(0.5,0.5){\nwlabel{$\lambda_F$}}
   \put(0.25,1){\vector(1,0){0.5}}
   \put(0.5,1){\rlabel{$\rho$}}
   \put(1,0){\object{$\moduli{r_1}{d_1}$}}
   \put(1.25,0){\vector(1,0){0.5}}
   \put(1.5,0){\rlabel{$\mu_1$}}
   \put(2,0){\object{$\moduli{h_1}{0}$}}
   \put(1,1){\object{$\pullback$}}
   \put(1.25,1){\vector(1,0){0.45}}
   \put(1,0.5){\dlabel{$\thetahat_1$}}
   \put(1,0.75){\vector(0,-1){0.5}}
   \put(1.5,1){\rlabel{$\muhat_1$}}
   \put(2,1){\object{$\parmod{h_1}{0}{h}$}}
   \put(2.3,1){\vector(1,0){0.4}}
   \put(2.5,1){\rlabel{$\theta_2$}}
   \put(2,0.75){\vector(0,-1){0.5}}
   \put(2,0.5){\dlabel{$\theta_1$}}
   \put(3,1){\object{$\moduli{h_1}{-h}$}}
   \put(3.3,1){\vector(1,0){0.45}}
   \put(3.5,1){\rlabel{$\mu_2$}}
   \put(4,1){\object{$\moduli{h}{0}$}}
  \end{picture}
\end{equation}
The maps $\mu_1$ and $\mu_2$ are of the same sort as $\mu$
and may be assumed to exist by induction, since both $r_1$ and 
$h_1=\hcf{r_1}{d_1}$ are less than $r$.
Thus they are birationally linear and satisfy
\begin{eqnarray}
  \mu_1^*(\brauer{h_1}{0})&=&\brauer{r_1}{d_1} \label{eq:mu1}\\
  \mu_2^*(\brauer{h}{0})&=&\brauer{h_1}{-h} \label{eq:mu2}.
\end{eqnarray}
because $h=\hcf{h_1}{-h}$ by Lemma~\ref{lem:r1d1}.

The central square in the diagram is a pull back.
In particular, $\widehat\theta_1\colon \pullback\ramto\moduli{r_1}{d_1}$
is the pull back of $\theta_1$ along $\mu_1$
and hence, by (\ref{eq:mu1}),
it is a Grassmannian bundle over $\moduli{r_1}{d_1}$
whose twisting is measured by $\brauer{r_1}{d_1}$.
Thus $\widehat\theta_1$ and $\lambda_F$ are twisted Grassmannian
bundles associated to vector bundles of weight $-1$ over $X_{r_1,d_1}$
and of ranks
$h_1$ and $lh_1$ respectively (see Proposition~\ref{t10}).
We will prove in Lemma~\ref{t4} below that this implies
that there is a birationally linear map
$\rho\colon\moduli{r}{d}\ramto\pullback$
such that $\lambda_F=\thetahat_1\rho$ and hence
\begin{equation}\label{eq:rho}
  \rho^*\left(\thetahat_1^*\left(\brauer{r_1}{d_1}\right)\right)
  =\brauer{r}{d}.  
\end{equation}
The pullback $\muhat_1$ of $\mu_1$ along $\theta_1$
is birationally linear and satisfies
\begin{equation}\label{eq:muhat1}
 \muhat_1^*\left(\theta_1^*(\brauer{h_1}{0})\right)
 =\thetahat_1^*\left(\mu_1^*(\brauer{h_1}{0}\right)
 =\thetahat_1^*\left(\brauer{r_1}{d_1}\right).
\end{equation}
Thus $\muhat_1\rho\colon \moduli{r}{d}\ramto \parmod{h_{1}}{0}{h}$ 
is birationally linear and pulls back $\theta_1^*(\brauer{h_1}{0})$
to $\brauer{r}{d}$.
But by (\ref{eq:thetas}) and (\ref{eq:mu2}), this means that 
$\mu=\mu_2\theta_2\muhat_1\rho$ pulls back $\brauer{h}{0}$
to $\brauer{r}{d}$ as required and to complete the proof
we need to show that $\theta_2$ is birationally linear.
This follows from Lemma~\ref{t3} below, because,
as we saw in Section~\ref{sec:Hecke}, 
$\theta_2$ is a twisted Grassmannian bundle of $h$-dimensional subspaces
of a vector bundle $H_2$ of weight 1 over $X_{h_1,-h}$
and the Brauer class $\brauer{h_1}{-h}$ is represented by a central
simple algebra of dimension $h^2$.
Thus although $\theta_2$ is not locally trivial in the Zariski
topology when $h\neq 1$, we can show that it is birationally linear
since its generic fibre is birational to a Grassmannian over a
division algebra; this is not the way it is expressed in
Lemma~\ref{t3} though the translation to this is fairly simple. 
\end{proof}

Thus (modulo two lemmas) we have proved Theorem~\ref{thm:1}
as we set out to do.
To deduce Theorem~\ref{thm:2}, it is sufficient to observe that,
by Lemma~\ref{lem:r1d1}, the map $\lambda_F$ restricts to a map between
moduli spaces of fixed determinant and that the Hecke correspondence
restricts to a correspondence between moduli spaces of fixed determinant.
Therefore the map $\mu\colon\moduli{r}{d}\to\moduli{h}{0}$
restricts to a map between fixed determinant moduli spaces,
although precisely how the determinants are related will depend
on various choices made in the construction.
In the case $h=1$, the fixed determinant moduli space is a point
and we obtain Theorem~\ref{thm:2}.

We finish with the proofs of the two lemmas about birationally linear maps
that we used in the proof of Theorem~\ref{thm:main}.
The first thing we need to understand is when a twisted Grassmannian bundle is
birationally linear over its base. We shall provide a sufficient
condition which is in fact necessary though we shall not prove that
here since we do not need it.

Start by observing that, if $P$ and $Q$ are vector bundles of weight
$w$ over $X$, then $M=P\dual\otimes Q/\PGL{n}$ is a bundle of
left modules for $P\dual\otimes P/\PGL{n}$ and that this correspondence
is invertible because $Q=P\otimes_{P\dual\otimes P}\gamma^{\ast}M$,
where $\gamma\colon X\to X/\PGL{n}$ is the quotient map.

\begin{lemma}
\label{t3}
Let $P$ be a vector bundle of weight $w$ over $X$.  Assume that the
Brauer class associated to $P$ is represented by a central simple
algebra of dimension $j^{2}$. Then $\pi\colon \Gr{j}{P}/\PGL{n}\rightarrow
X/\PGL{n}$ is a birationally linear map.
\end{lemma}
\begin{proof}
Let $A$ be the central simple algebra given by the bundle of central
simple algebras $P\dual\otimes P/\PGL{n}$ over the field
$\funfield{X/\PGL{n}}$. Then by assumption $A$ has a left ideal
of dimension $j\rk(P)$ which is of necessity a direct summand of
$A$. Therefore, over some dense open subvariety of $X/\PGL{n}$,
$P\dual\otimes P/\PGL{n} \cong L_{1}\oplus L_{2}$ where $L_{1}$ and
$L_{2}$ are bundles of left ideals for $P\dual\otimes P/\PGL{n}$ and
$\rk(L_{1})=j\rk(P)$. We may as well assume that this happens over
$X/\PGL{n}$. We obtain a corresponding direct sum decomposition of
$P$, $P \cong P_{1}\oplus P_{2}$ where $P_{1}$ and $P_{2}$ are
$\GL{n}$ stable subbundles of $P$ and hence both of weight $w$. Also
$\rk(P_{1})=j$. Now consider the vector bundle 
$P_{1}\dual\otimes P_{2}$. Let $\lambda\colon P_{1}\rightarrow P_{2}$ be the
universal homomorphism of vector bundles defined on 
$P_{1}\dual\otimes P_{2}$ and consider the map of vector bundles
over $P_{1}\dual\otimes P_{2}$, 
$(Id,\lambda)\colon P_{1}\rightarrow P_{1}\oplus P_{2}\cong P$. This
representation of $P_{1}$ as a subbundle of $P$ defines a map from 
$P_{1}\dual\otimes P_{2}$ to $\Gr{j}{P}$ which is
$\PGL{n}$-equivariant, injective and onto an open subvariety of
$\Gr{j}{P}$. Hence $P_{1}\dual\otimes P_{2}/\PGL{n}$ which is a vector
bundle over $X/\PGL{n}$ is an open subvariety of $\Gr{j}{P}$. 
\end{proof}

It remains to show that two twisted Grassmannian bundles of
equal dimensional subspaces arising from
vector bundles of the same weight
have a birationally linear map between them.
 
\begin{lemma}
\label{t4}
Let $P$ and $Q$ be vector bundles of weight $w$ over $X$
and suppose that $j<\rk(Q)<\rk(P)$.
Then there is a birationally linear rational map
$$\rho\colon\Gr{j}{P}/\PGL{n}\to\Gr{j}{Q}/\PGL{n}.$$
compatible with the bundle maps to $X/\PGL{n}$.
\end{lemma}
\begin{proof}
$P\dual\otimes Q/\PGL{n}$ is a bundle of left modules for
$P\dual\otimes P/\PGL{n}$ of rank equal to
$\rk(P)\rk(Q)$ and since
$\rk(Q)<\rk(P)$ there is an open subvariety of $X/\PGL{n}$ on which
$$P\dual\otimes P/\PGL{n} \cong P\dual\otimes Q/\PGL{n}\oplus L$$
for some vector bundle of left ideals $L$ since this is true at the
generic point of $X/\PGL{n}$. Hence we may assume that on
$X/\PGL{n}$, $P\cong Q\oplus Q\dsh$ for $\GL{n}$ stable subbundles $Q$
and $Q\dsh$. Let $S$ be the universal subbundle on $\Gr{j}{Q}$. Then
$S$ and $Q\dsh$ are both vector bundles of weight $w$ on
$\Gr{j}{Q}$. We consider the vector bundle $S\dual\otimes Q\dsh$ over
$\Gr{j}{Q}$. Let $\lambda\colon S\rightarrow Q\dsh$ be the universal
homomorphism of vector bundles defined on $S\dual\otimes Q\dsh$ and
let $\iota\colon S\rightarrow Q$ be the universal inclusion of $S$ in $Q$
pulled back to $S\dual\otimes Q\dsh$; now consider the map of
vector bundles 
$$(\iota,\lambda)\colon S\rightarrow Q\oplus Q\dsh\cong P$$
defined on $S\dual\otimes Q\dsh$. This gives a subbundle of $P$ of rank
$j$ and hence defines a map from $S\dual\otimes Q'$ to
$\Gr{j}{P}$. This map is injective and onto an open subvariety of
$\Gr{j}{P}$ and it is also $\PGL{n}$-equivariant. Hence 
$S\dual\otimes Q\dsh/\PGL{n}$ is an open subvariety of
$\Gr{j}{P}/\PGL{n}$. However, $S\dual\otimes Q\dsh/\PGL{n}$ is a
vector bundle over $\Gr{j}{Q}/\PGL{n}$ which gives our lemma. 
\end{proof}


\end{document}